\documentclass{alf-j-l}
\usepackage{verbatim} 
\usepackage{syntonly}
\usepackage{url}  

\def\Z{{\mathbb Z}}  
\def\F{{\mathbb F}} 

\def\\{\cr} 
\def\lt{<}
\def\gt{>}
\def\conj{\overline}

\def\notdivides{\mathrel{\kern-3pt\not\!\kern3.5pt\bigm|}}

\def\poly {polynomial}
\def\cf {continued fraction}
\def\pq {partial quotient}
\def\cq {complete quotient}

\def\ex{expansion}

\def\cfe{continued fraction expansion}
\DeclareMathOperator{\Jac}{Jac}

\theoremstyle{plain}  
\newtheorem{theorem}{Theorem}

\newtheorem*{theorem*}{Theorem}

\theoremstyle{definition}

\newtheorem*{Remark*}{Remark}
\newtheorem*{Remarks*}{Remarks}

\newtheorem*{example*}{Example} 
\newtheorem*{guess*}{Guess}

\theoremstyle{remark}

\newtheorem*{remark*}{Remark} 
\newtheorem*{remarks*}{Remarks}

\newtheoremstyle{aside}
   {6pt}
   {6pt}
   {\footnotesize}
   {}
   {\scshape}
   {:}
   {.5em}
   {}

\theoremstyle{aside}

\newtheorem*{aside*}{Concluding Aside}

\begin{document}

\title {Curves of Genus~2,
Continued Fractions, \\and Somos Sequences} 

\author{Alfred J. van der Poorten}
\address{Centre for Number Theory Research, 1 Bimbil Place, Killara, Sydney,
NSW 2071, Australia}
\email{alf@math.mq.edu.au (Alf van der Poorten AM)}

\thanks{The author's only support was a grant from the Australian
Research Council.}



\begin{abstract} We detail the continued fraction expansion of the square
root of monic sextic polynomials. We note in passing that each line
of the expansion corresponds to addition of the divisor at infinity,
and interpret the data yielded by the general expansion. In
particular we obtain an associated Somos sequence defined by a
three-term recurrence relation of gap~$6$.
\end{abstract}

\maketitle

\noindent In the present note I study the \cfe\ of the square root of a sextic
\poly, inter alia obtaining integer sequences generated by recursions
$$A_{h-3}A_{h+3}=aA_{h-2}A_{h+2}+bA_{h}^2\,.
$$
Specifically, see \S\ref{s:example} at page~\pageref{s:example} for the case 
$(T_h)=(\ldots, 2,1,1,1,1,1, 1, 2,3,\ldots)$,  where I illustrate how the \cfe\ data
readily allows one to recover the genus~$2$ curve $\mathcal C:
Y^2=(X^3-4X+1)^2+4(X-2)$ giving rise to the sequence. 

\section{Some Brief Reminders}
\noindent A reminder exposition on \cf s in quadratic function fields appears as
\S4 of \cite{169}. However, the na\"{\i}ve reader needs little more than that a \cfe\ of a
quadratic irrational integer $\omega$ is a two-sided sequence of lines, $h$ in $\Z$,
\begin{equation*}
\frac{\omega+P_h}{Q_h}=a_h-\frac{\conj\omega+P_{h+1}}{Q_h}\,;\quad\text{in brief
$\omega_h=a_h-\conj \rho_h$}\,,
\end{equation*}
with $(\omega+P_{h+1})(\conj\omega+P_{h+1})=-Q_hQ_{h+1}$ defining the integer sequences $(P_h)$
and $(Q_h)$. Necessarily, one must have, say, $Q_0$ divides
$(\omega+P_{0})(\conj\omega+P_{0})$ in which case the integrality of the sequence $(a_h)$ of
\emph{\pq s} guarantees that always $Q_h$ divides $(\omega+P_{h})(\conj\omega+P_{h})$. 
If the \emph{\pq} $a_h$ is always chosen as the integer part of
$\omega_h$ then
$\omega_0$ reduced implies that all the
$\omega_h$ and $\rho_h$ are reduced; and always $a_h$ also is the integer part of
$\rho_h$. Then conjugation retrieves the left hand half of the \ex\ of $\omega_0$ from
that of $\rho_0$. In the function field case, one reads `\poly' for `integer'.
\section{Continued Fraction Expansion of the Square Root of a Sextic}
\label{s:cfe}

\noindent We suppose the base field $\F$ is not of characteristic~$2$ or $3$
because those cases requires changes throughout the exposition; indeed, nontrivial
changes in the case characteristic~2. We study the \cfe\ of the
squre root of a sextic \poly\
$D\in\F[X]$. Set
\begin{equation} \label{eq:C}
\mathcal C: Y^2=D(X):=(X^3+fX+g)^2+4u(X^2-vX+w),
\end{equation}
and for brevity write $A=X^3+fX+g$ and $R=u(X^2-vX+w)$. Set
$Z=\frac12(Y+A)$ and notice that
$Z^2-AZ-R=0$. The other root of this equation is $\conj Z$. 

Suppose that $(X^2-v_0X+w_0)$ divides the norm
$$Z_0\conj Z_0=-R+d_0(X+e_0)\bigl(A++d_0(X+e_0)\bigr)\,,$$
and that $Z_0$ has been so chosen that all\footnote{At any rate,
sufficiently many \pq s so as to make the following discussion
useful.} its \pq s are of degree~$1$. Such a choice is `generic' if
the base field is infinite.

For 
$h=0$,
$1$,
$2$,
$\ldots\,$  we denote the \cq s of $Z_0$ by
\begin{equation}
\label{eq:notation}Z_h=\bigl(Z+d_h(X+e_h)\bigr)/u_h(X^2-v_hX+w_h)\,,
\end{equation}
noting that the $Z_h$ all are reduced, namely $\deg Z_h\gt0$ but $\deg
\conj Z_h\lt0$. The upshot is that the $h$-th line of the
\cfe\ of
$Z_0$ is
\begin{equation}\label{eq:P}
\frac{Z+d_h(X+e_h)}{u_h(X^2-v_hX+w_h)}=\frac{X+v_h}{u_h} - 
\frac{\conj Z+d_{h+1}(X+e_{h+1})}{u_h(X^2-v_hX+w_h)}\,.
\end{equation}
Then evident recursion formulas yield
\begin{equation} \label{eq:d} 
f+d_h+d_{h+1}=-v_h^2+w_h
\end{equation}
\begin{equation} \label{eq:e} 
g+d_he_h+d_{h+1}e_{h+1}=v_hw_h
\end{equation}
and 
\begin{multline}\label{eq:Q}
-u_hu_{h+1}(X^2-v_hX+w_h)(X^2-v_{h+1}X+w_{h+1})\\
=\bigl(Z+d_{h+1}(X+e_{h+1})\bigr)\bigl(\conj Z+d_{h+1}(X+e_{h+1})\bigr).
\end{multline}
Hence, noting that $Z\conj Z=-u(X^2-vX+w)$ and $Z+\conj Z=A=X^3+fX+g$, we may
equate coefficients in \eqref{eq:Q} to see that
\begin{equation} \label{eq:X4}
d_{h+1}=-u_hu_{h+1}\,.
\tag{$\ref{eq:Q}:X^4$} \end{equation}
Given that, we obtain, after in each case dividing by
$-u_hu_{h+1}$,
\begin{equation} \label{eq:X3}
e_{h+1}=-v_h-v_{h+1}\,;
\tag{$\ref{eq:Q}:X^3$} \end{equation}
\begin{equation} \label{eq:X2}
(f+d_{h+1})=v_hv_{h+1}+(w_h+w_{h+1})+u/d_{h+1}\,;
\tag{$\ref{eq:Q}:X^2$} \end{equation}
\begin{equation} \label{eq:X1}
(f+d_{h+1})e_{h+1}+(g+d_{h+1}e_{h+1})
=-v_hw_{h+1}-v_{h+1}w_h-uv/d_{h+1}\,;
\tag{$\ref{eq:Q}:X^1$} \end{equation}
\begin{equation} \label{eq:X0}
(g+d_{h+1}e_{h+1})e_{h+1}
=w_hw_{h+1}+uw/d_{h+1}\,.
\tag{$\ref{eq:Q}:X^0$} \end{equation}

The :$X^2$ equation readily becomes
\begin{equation*}
-d_h=f-w_h+v_h^2+d_{h+1}=v_h(v_h+v_{h+1})+w_{h+1}+u/d_{h+1}\,,
\end{equation*}
so $d_{h+1}(v_he_{h+1}-w_{h+1})=d_hd_{h+1}+u$. With similar manipulation of the
next two equations we felicitously obtain
\begin{subequations}\label{eq:helpful}
\begin{align} 
d_{h+1}(v_he_{h+1}-w_{h+1})&=d_hd_{h+1}+u\,;\\
-v_hd_{h+1}(v_he_{h+1}-w_{h+1})&=d_hd_{h+1}(e_h+e_{h+1})-uv\,;\\
w_{h}d_{h+1}(v_he_{h+1}-w_{h+1})&=d_hd_{h+1}e_he_{h+1}+uw\,.
\end{align}
\end{subequations}
That immediately yields
\begin{subequations}\label{eq:useful}
\begin{align} 
d_hd_{h+1}(e_h+e_{h+1}+v_h)&=u(v-v_h)\,;\\
d_hd_{h+1}(e_he_{h+1}-w_h)&=-u(w-w_h)\,.
\end{align}
\end{subequations}
Incidentally, by 
\begin{equation*}
-d_{h+1}=f-w_h+v_h^2+d_{h}=v_h(v_{h-1}+v_{h})+w_{h-1}+u/d_{h}\,,
\end{equation*}
we also discover that, mildly surprisingly,
\begin{equation}\label{eq:surprise}
d_hd_{h+1}+u=d_{h+1}(v_he_{h+1}-w_{h+1})=d_{h}(v_he_{h}-w_{h-1}).
\end{equation}

\section{A Ridiculous Computation}\label{s:ridiculous}

\noindent It is straightforward to notice that the three final equations
\eqref{eq:Q} yield
$$
e_h^2(v_{h-1}v_h+w_{h-1}+w_h)+e_h(v_{h-1}w_h+v_hw_{h-1})+w_{h-1}w_h=
-u(e_h^2+ve_h+w)/d_h\,.
$$
Remarkably, by \eqref{eq:surprise}
\begin{multline*}
(d_{h-1}d_h+u)(d_hd_{h+1}+u)=d_h^2(v_{h-1}e_{h}-w_{h})(v_he_{h}-w_{h-1})\\
=e_h^2v_{h-1}v_h-e_h(v_{h-1}w_{h-1}+v_hw_h)+w_{h-1}w_h
\end{multline*}
and so, because
\begin{multline*}
-(v_{h-1}w_{h-1}+v_hw_h)=v_{h-1}w_h+v_hw_{h-1}-(w_{h-1}+w_h)(v_{h-1}+v_h)\\
=v_{h-1}w_h+v_hw_{h-1}+e_h(w_{h-1}+w_h)\,,
\end{multline*}
we obtain the surely useful identity
\begin{equation}\label{eq:identity}
(d_{h-1}d_h+u)(d_hd_{h+1}+u)=-ud_h(e_h^2+ve_h+w)\,.
\end{equation}
This just one of the nine such identities provided by the equations
\eqref{eq:helpful}, and
\eqref{eq:surprise}.

\subsection{The special case $u=0$} Consider now the case in which $R$, the
remainder term
$u(X^2-vX+w)$, is replaced by $-v(X-w)$. In effect $u\gets0$ except that $uv\gets
v$, $uw\gets vw$. For instance, \eqref{eq:identity} becomes
 \begin{equation}\label{eq:identity'}
d_{h-1}d_hd_{h+1}=-v(e_h+w)\,,\tag{\ref{eq:identity}$'$}
\end{equation}
and, we'll need this, we now have
\begin{align} 
\label{eq:v} e_h+e_{h+1}+v_h&=v/d_hd_{h+1}\,;\tag{\ref{eq:useful}$'$a}\\
\label{eq:w} e_he_{h+1}-w_h&=-vw/d_hd_{h+1}\,.\tag{\ref{eq:useful}$'$b}
\end{align}
Indeed, we find that
\begin{equation}
d_{h-1}d_h^2d_{h+1}^2d_{h+2}=v^2(e_he_{h+1}+w(e_h+e_{h+1})+w^2)=
v^2(w_h-wv_h+w^2)
\end{equation}
and therefore that
\begin{multline}
d_{h-2}d_{h-1}^3d_h^4d_{h+1}^3d_{h+2}=\\
v^4\bigl(w_{h-1}w_h+w^2\bigl(v_{h-1}v_h+(w_{h-1}+w_h)\bigr)-w(v_{h-1}w_h+w_{h-1}v_{h})-
w^3(v_{h-1}+v_h)+w^4\bigr)\,.
\end{multline}
This last expression is transformed by the equations \eqref{eq:Q} to become
\begin{multline}
v^4\bigl((g+d_he_h)e_h-vw/d_h+w^2(f+d_h)+
\\+w\bigl((f+d_h)e_h+(g+d_he_h)+v/d_h\bigr)+w^3e_h+w^4\bigr)\\
=v^4(e_h+w)\bigl((g+d_he_h)+w(f+d_h)+w^3\bigr)\,.
\end{multline}
Thus
\begin{equation}
d_{h-2}d_{h-1}^2d_h^3d_{h+1}^2d_{h+2}=-v^3\bigl((g+d_he_h)+w(f+d_h)+w^3\bigr)\,.
\end{equation}
But wait, there's more! By \eqref{eq:identity'} we know that $-ve_h=d_{h-1}d_hd_{h+1}+vw$, so
\begin{equation}\label{eq:more}
d_{h-2}d_{h-1}^2d_h^3d_{h+1}^2d_{h+2}=v^2d_{h-1}d_h^2d_{h+1}-v^3(g+wf+w^3)\,.
\end{equation}

\begin{theorem} Set $D(X)=A^2+4R:=(X^3+fX+g)^2-4v(X-w)$ and let $Z=\frac12(Y+A)$, so
$Z^2-AZ-R=0$. Denote by
$$Z_h=\frac{Z+d_h(X+e_h)}{u_h(X^2-v_hX+w_h)}\,,
$$
$h\in\Z$, the complete quotients of the \cfe\ of $Z_0$; here $Q_0(X)=u_0(X^2-v_0X+w_0)$ must
divide the norm $d_0^2(X+e_0)^2+d_0(X+e_0)A-R\,.$
\end{theorem}
Denote by $(T_h)$ a sequence defined by appropriate initial values and the recursive relation
\begin{equation}\label{eq:relation}
T_{h-1}T_{h+1}=d_hT_h^2\,.
\end{equation}
Then
\begin{equation}\label{eq:recurrence}
T_{h-3}T_{h+3}=v^2T_{h-2}T_{h+2}-v^3(g+wf+w^3)T_h^2\,.
\end{equation}
\begin{proof}
It suffices to check that, given \eqref{eq:relation}, we need only
multiply \eqref{eq:more} by $T_h^2$.
\end{proof}
\begin{Remark*} The reader should note the evident tight analogy with the
corresponding result for quartic \poly s detailed in \cite{169}.
On the other hand, the results of \cite{169} continue to
make sense even in \emph{singular} cases,
when there are \pq s of degree greater than~one. That's
not quite so here: surely, $T_{k-3}=T_{k-2}=0$ is usually not
compatible with \eqref{eq:recurrence}, suggesting that
then our division by $e_k+w$ in the course of our
`ridiculous argument' may be an improper division by
zero.\end{Remark*}

\section{A Cute Example}\label{s:example}

\noindent The example
\begin{equation}\label{eq:example}
T_{h-3}T_{h+3}=T_{h-2}T_{h+2}+T_h^2\,,
\end{equation}
with $T_0=T_1=T_2=T_3=T_4=T_5=1$ is readily seen to derive from the genus~$2$ curve
\begin{equation}\label{eq:examplecurve}
\mathcal C: Y^2=(X^3-4X+1)^2+4(X-2)\,.
\end{equation}
To indeed see this, we first note that of course we need $d_1=d_2=d_3=d_4=1$ to produce the
initial values from
$T_0=T_1=1$. Since, plainly, $T_{-1}=T_6=2$, clearly $d_0=2$. By the Theorem, we expect to
require
$v^2=1$ and
$-v^3(g+wf+w^3)=1$. Without loss of generality, we may take $v=-1$. From
\eqref{eq:identity'} we then read off that
\begin{equation*}
e_1=2-w \quad\text{and}\quad e_2=1-w\,.
\end{equation*}
Thus, by \eqref{eq:d} and \eqref{eq:e} we have
$$
f+2=-v_1^2+w_1 \quad\text{and}\quad g +3-2w=v_1w_1\,.
$$
But from \eqref{eq:v} and \eqref{eq:w} we evaluate $v_1$ and $w_1$ in terms of $w$ as 
$$3-2w+v_1=-1 \quad\text{and}\quad (2-w)(1-w)-w_1= w\,.
$$
Substituting appropriately we find that $1=g+fw+w^3=6w-11$ so, as already announced,
$v=-1$, $w=2$,
$g=1$, and
$f=-4$. 

Furthermore, we have $v_1=0$ and $v_0+v_1+e_1=0$, so $v_0=0$; then $f+3=-v_0^2+w_0$ yields
$w_0=-1$. Noting that $g+2e_0+e_1=0$, we find that $e_0=-1/2$. Thus the relevant \cfe\ commences
\begin{align*}
Z_0:=\frac{Z+2X-1}{X^2-1}&=X-\frac{\conj Z+X}{X^2-1}\\
\frac{Z+X}{-(X^2-2)}&=-X-\frac{\conj Z+X-1}{-(X^2-2)}\\
\frac{Z+X-1}{X^2-X-1}&=X+1-\frac{\conj Z+X-1}{X^2-X-1}\\
\frac{Z+X-1}{-(X^2-2)}&=-X-\frac{\conj Z+X}{-(X^2-2)}\\
\frac{Z+X}{X^2-1}&=X-\frac{\conj Z+2X-1}{X^2-1}\\
&\cdots\end{align*}
providing a useful check on our allegations and displaying an expected symmetry.
Denote by $M$ the divisor class defined by the pair of points
$(\varphi,0)$ and
$(\conj\varphi,0)$ --- here, $\varphi$ is the golden ratio, a
happenstance that will please adherents to the cult of Fibonacci ---
and by
$S$  the divisor class at infinity. Then the sequence
$(T_h)=(\ldots, 2, 1,1,1,1,1, 1, 2,3,4,8,17,50
\ldots)$ may be thought of as arising from the points $\ldots$,
$M-S$, $M$, $M+S$, $M+2S$, $\ldots$ on the Jacobian of the curve
$\mathcal C$ displayed at \eqref{eq:examplecurve}. Evidently,
$M-S=-M$ so $2M=S$ on $\Jac(\mathcal C)$. 

Incidentally, this closing aside identifying the \cfe\ with stepwise addition on the
Jacobian is gratuitous. However, concerned readers might contemplate the introduction to Cantor's
paper
\cite{Ca} and the instructive discussion by Kristin Lauter in \cite{La}. A central theme of
the paper \cite{BCZ} is a generalisation of the phenomenon to Pad\'e approximation in
arbitrary algebraic function fields.

\section{Comments}

\noindent I consider the argument given in \S\ref{s:ridiculous}  above to be quite absurd
and am ashamed to have spent a great deal of time in extracting it. Such are the costs of
truly low lowbrow arguments; see \cite{CF} for heights of `brow'. The only saving grace is
my mildly ingenious use of symmetry in the argument's later stages. I do not know whether
there is an appealing result of the present genre if
$u\ne0$; but see my remarks below. I should admit that I realised, but only after having
successfully selected
$u=0$, that Noam Elkies had suggested to me at ANTS, Sydney 2002, that an identity of the
genre
\eqref{eq:recurrence} would exist, but had in fact specified just the special case $\deg
R=1$.

Mind you, with some uninteresting effort one can show (say by counting free
parameters) that over an algebraic extension of the base field there is a birational
transformation which transforms the given curve to one where $\deg R=1$. That does not
truly better the present theorem.

On the other hand, a dozen years ago\footnote{I have a revision of his manuscript dated
November, 1992.}, David Cantor \cite{Ca} mentions that his results lead readily to Somos
sequences both in genus~$1$ and $2$; the latter with gap~$8$ (see \cite{169} for the
relevant notions). The latter consequence is not obvious; however, recently, Cantor has
told me a rather ingenious idea which clearly yields the result for all
hyperelliptic curves $Y^2=E(X)$, $E$ a quintic, say with constant coefficient~$1$. In
brief, Cantor's result is more general than mine but does not deal with all cases I handle
here; nor does it produce the expected recursion formul\ae\ of gap~$6$.

The most serious disappointment is that the best argument I can produce here is just a
much more complex version of that of \cite{169} for genus~$1$. Seemingly, a different
\emph{Ansatz}, a new view on the issues, is needed if my methods are
to yield results in higher genus.

\bibliographystyle{amsalpha}

\label{page:lastpage}
\end{document}